\documentclass [12pt]{article}
\usepackage {amsfonts}
\usepackage {mathrsfs}
\usepackage {amsmath}
\usepackage {latexsym}
\usepackage {amssymb}
\usepackage {amsthm}
\usepackage {amscd}
\date{}
\title{\bf the relationship between  graphs
and Nichols braided Lie algebras}
\author{\small Weicai Wu $^{a}$, Shouchuan Zhang $^{b}$, Zhengtang Tan$^{c}$\\
\small $a$.  School of Mathematics, Hunan Institute of Science and Technology,\\
\small Yueyang  414006, P.
 R. China.\\
\small $b$. Department  of Mathematics,   Hunan University,\\
\small   Changsha  410082,   \ P.R. China.\\
\small $^c$. School of Engineering and Design, Hunan Normal University,\\
\small Changsha  $410081$,   P.R. China.\\
\small \tt Emails: z9491@sina.cn;  weicaiwu@hnu.edu.cn  }

\begin{document}
\newtheorem{Proposition}{Proposition}[section]
\newtheorem{Theorem}[Proposition]{Theorem}
\newtheorem{Definition}[Proposition]{Definition}
\newtheorem{Corollary}[Proposition]{Corollary}
\newtheorem{Lemma}[Proposition]{Lemma}
\newtheorem{Example}[Proposition]{Example}
\newtheorem{Remark}[Proposition]{Remark}

\maketitle %\addtocounter{section}{-1}

\begin {abstract}

In this paper we give the relationship between the connected components of  pure generalized Dynkin graphs
and Nichols braided Lie algebras.


\end {abstract}

\section {Introduction}\label {s0}
Nichols algebras have found significant applications in various areas of mathematics and mathematical physics including the
theories of pointed Hopf algebras and logarithmic quantum fields. In \cite {He05}, one-to-one correspondences
between Nichols algebras of diagonal type and arithmetic root systems as well as between generalized Dynkin diagrams
and twisted equivalence classes of arithmetic root systems were established.
The problem of finite-dimensionality of Nichols algebras forms a substantial part of the recent investigations
(see e.g. \cite {AHS08,  AS10,  He05,  He06a,  He06b,  WZZ15a,  WZZ15b}. Lie algebra arising from a Nichols algebra was studied in \cite {AAB16}.

Braided Lie algebras  were studied in \cite {Ka77, BMZP92, Gu86, GRR95, Kh99, Pa98, Sc79, BFM96, ZZ03}. The current paper will focus on
Nichols  Lie and braided Lie algebras.  In \cite {He05} and \cite {He06a}, a classification on braided vector spaces of diagonal type
with finite-dimensional Nichols algebras was given. In \cite {WZZ15b, WZZ16}, we studied the relationship between Nichols algebras and Nichols braided Lie algebras.
It was proven that a Nichols algebra is finite-dimensional if and only if the corresponding Nichols (braided) Lie algebra is finite-dimensional. In \cite {WWZZ17} we show that a monomial belongs to Nichols braided Lie algebra $\mathfrak L(V)$ of braided vector space $V$ of diagonal type
if and only if that this monomial is connected. This is one of the main results in this paper, which enables us to obtain the bases of the
Nichols braided Lie algebras of arithmetic root systems and the dimensions of Nichols braided Lie algebras of finite Cartan type.
We give the sufficient and necessary conditions for  $\mathfrak B(V) = F\oplus \mathfrak L^-(V)$ and $\mathfrak L^-(V)= \mathfrak L(V)$, where
$\mathfrak B(V)$, $\mathfrak L^-(V)$ and $\mathfrak L(V)$ denote Nichols algebra, Nichols Lie algebra
and Nichols braided Lie algebra over $V$. We also obtain an explicit basis of $\mathfrak L^ - (V)$ over quantum linear space $V$ with $\dim V=2$.
This provides a new method for determining when a Nichols algebra is finite dimensional.

In this paper we give the relationship between the connected components of  pure generalized Dynkin graphs
and Nichols braided Lie algebras;

\section*{ Preliminaries}

For any matrix $(q_{ij})_{n \times n}$  over $F^*$, define a bicharacter $\chi$ from $\mathbb Z^n \otimes \mathbb Z^n$ to $F^*$ such that $\chi (e_i, e_j) =q_{ij}$  for $1\le i, j \le n$, where $\{ e_1, e_2, \cdots, e_n\}$ is a basis of $\mathbb Z^n.
$ Let $V$ be a vector space  with basis
$x_1,  x_2,  \cdots,  x_n$. Define $\alpha (e_i, x_j) = q_{ij}x_j$ and $\delta (x_j) =e_j \otimes x_j$ for $1\le i, j \le n$. It is clear that $(V, \alpha, \delta)$ is a Yetter-Drinfeld module over $\mathbb Z^n$ and $(V, C) $ is a  braided vector space under braiding $C,$  where $C(x_i \otimes x_j) = q_{ij}x_j \otimes x_i$
and $C^{-1}(x_i \otimes x_j) = q_{ji} ^{-1}x_j \otimes x_i.$ In this case, $V$ is called a  braided vector space  of diagonal type and $(q_{ij})_{n \times n}$ is called a braiding matrix of $V$.
Throughout this paper  braided vector space $V$ is of diagonal type with basis
$x_1,  x_2,  \cdots,  x_n$ and $C(x_i \otimes x_j) =q_{ij} x_j \otimes x_i$ without special announcement. Let $\mathfrak B(V)$ be the Nichols algebra over the braided  vector space $V$.
Define $p_{ij} := q_{ij}$ for $1\le i, j \le n$  and  $p_{uv} := \chi ( {\rm deg} (u),   {\rm deg }(v))  $ for any homogeneous element $u,   v \in \mathfrak B(V).$
Denote ${\rm ord } (p_{uu})$ the  order of $p_{uu}$ with respect to multiplication. Let $|u|$ denote length of  homogeneous element $u\in \mathfrak B(V).$
Let $D =: \{[u] \mid [u] \hbox { is a hard super-letter }\}$, $\Delta ^+(\mathfrak B(V)): =  \{ \deg (u) \mid [u]\in D\}$,
$ \Delta (\mathfrak B(V)) := \Delta ^+(\mathfrak B(V)) \cup \Delta ^-(\mathfrak B(V))$,
which is called the root system of $V.$
 If $ \Delta (\mathfrak B(V))$ is finite,
then it is called an arithmetic root system. Let  $\mathfrak L(V)$ denote the  braided Lie algebras generated by $V$ in $\mathfrak B(V)$ under Lie operations $[x, y]=yx-p_{yx}xy$,  for any homogeneous elements $x, y\in \mathfrak B(V)$. $(\mathfrak L(V), [\ ])$ is called Nichols braided Lie algebra of $V$. Let  $\mathfrak L^-(V)$ denote the  Lie algebras generated by $V$ in $\mathfrak B(V)$ under Lie operations $[x, y]^-=yx-xy$,  for any homogeneous elements $x, y\in \mathfrak B(V)$. $(\mathfrak L^-(V), [\ ]^-)$ is called Nichols  Lie algebra of $V$.
The other notations are the same as in \cite {WZZ15a}.

\vskip.1in
Recall the dual $\mathfrak B(V^*) $ of Nichols algebra $\mathfrak B(V) $ of rank $n$
in \cite [Section 1.3]{He05}. Let $y_{i}$ be a dual basis of $x_{i}$. $\delta (y_i)=g_i^{-1} \otimes y_i$,  $g_i\cdot y_j=p_{ij}^{-1}y_j$ and $\Delta (y_i)=g_i^{-1}\otimes y_i+y_i\otimes 1.$ There exists a bilinear map
$<\cdot, \cdot >$: $(\mathfrak B(V^*)\# FG)\times\mathfrak B(V)$ $\longrightarrow$ $\mathfrak B(V)$ such that

$<y_i, uv>=<y_i, u>v +g_i^{-1}.u<y_i, v>$ and $<y_i, <y_j, u>>=<y_iy_j, u>$

\noindent for any $u, v\in\mathfrak B(V)$. Furthermore,  for any $u\in \oplus _{i=1}^\infty \mathfrak B(V)_{(i)}$,
one has  $u=0$ if and only if $<y_i, u>=0$ for any $1\leq i \leq n.$
We have
\begin {eqnarray}\label {e0.1} [[u, v], w]=[u, [v, w]]+p_{vw}^{ -1} [[u, w], v]+(p_{wv}-p_{vw}^{-1})v\cdot[u, w],
\end {eqnarray}
\begin {eqnarray}\label {e0.2} [u, v \cdot w]=p_{wu}[uv]\cdot w + v \cdot [uw].\end {eqnarray}

We now recall some basic concepts of the graph theory (see \cite {Ha69}).
Let $\Gamma _1$ be a non-empty set and $\Gamma _2 \subseteq \{  \{ u, v\} \mid u, v \in \Gamma _1, \hbox { with } u \not= v \} \subseteq 2 ^{\Gamma_1}.$ Then $\Gamma = (\Gamma _1, \Gamma_2)$ is called a graph; $\Gamma_1$ is called the vertex set of  $\Gamma$; $\Gamma_2$ is called the edge set of  $\Gamma$; Element $\{u, v\} \in \Gamma_2$ is called an edge, written $a_{u, v}$.
 If $G = (G_1, G_2)$ is a graph and $G_1 \subseteq \Gamma _1$ and $G_2 \subseteq \Gamma _2$, then $G$ is called a subgraph of $\Gamma.$ If $ \emptyset \not= H_1 \subseteq \Gamma_1$ and $H_2 = \{ a_{u, v} \in \Gamma _2\mid u, v \in H_1 \}$, then $H=(H_1, H_2)$ is a subgraph, called the subgraph generated by $H_1$ in $\Gamma$.

 $a_{u_mu_{m-1}} \cdots a_{u_3u_2} a_{u_2u_1}$ is called a walk  from $u_1$ to $u_m.$
 We can define an equivalent relation on $\Gamma_1$ as follows: for any $u, v\in \Gamma _1$, $u$ and $v$ are equivalent if and only if there exists a walk from $u$ to $v$ or $u=v.$  Every subgraph  generated by  every equivalent class of $\Gamma_1$ is called
a connected component of $\Gamma.$

Let $\Gamma (V)$ be the  generalized Dynkin diagram of $V$ which excludes $p_{x_i,  x_i},\widetilde{p}_{x_i,  x_i}$ for $1\le i \le n $. This is called a pure generalized Dynkin graph of $V$, i.e. $\Gamma (V)_1 = \{ 1, 2, \cdots, n\}$ and $\Gamma (V)_2 = \{  a_{ij} \mid  p_{ij} p_{ji} \not= 1, i \not= j\}$.

Let
$\Gamma_a (V)_1 = \{ 1, 2, \cdots, n\}$ and $\Gamma_a (V)_2 = \{  a_{ij} \mid  p_{ij} \not= 1, i \not= j\}$. $\Gamma_a (V) = (\Gamma_a (V)_1, \Gamma_a (V)_2)$ is called the augmented Dynkin graph of $V$.

Let   $u=  h_1 h_2\cdots  h_m$ be   a  monomial with $ h_j = x_{i_j}$  for $1\leq i_1,i_2,\ldots,i_m\leq n,1\le j \le m$. ${\rm deg} (u) = \lambda _1 e_1 + \cdots + \lambda _n e_n, $ where $ {\rm deg }(x_i) = e_i.$ Let ${\rm deg}_{x_i} (u) := \lambda _i$. Let $\mu (u):= \{ x_{i_1},  \cdots,  x_{i_m}\}$  and  $\Gamma (u)$ be a pure generalized Dynkin subgraph generated by $\mu (u)$.
 If $\Gamma (u)$
is connected,  then $u$ is called connected (or $\mu (u) $ is called connected). Otherwise, $u$ is called disconnected (or $\mu (u) $ is called disconnected).

For $u, v \in \mathfrak B(V)$, if there exists a non-zero $a\in F$ such that  $u = av$, then we write  $u \sim v.$ This is an equivalent relation.
If there exist $x_i \in \mu (u)$ and   $x_j \in \mu (v)$   such that $\widetilde{p} _{x_i,  x_j} \not=1$,  then we say that it is connected between monomial $u$ and monomial $v$, written $u \diamondsuit v$ in short. Otherwise, we say that it is disconnected between monomial $u$ and monomial $v$.

Remark:
When $u\not=0$, $\mu (u)$ is independent of the choice of $h_1, h_2, \cdots, h_m,$ since $ \mathfrak B(V)$ is graded and ${\rm deg } (u) $ is unique. Therefore, the connectivity of $u$ is independent of the choice of $h_1, h_2, \cdots, h_m.$ When $u=0$, $\mu (u)$ is dependent on the choice of $h_1, h_2, \cdots, h_m.$  For example, ${\rm deg } (V) >1$, $p_{12} p_{21} =1, $ $p_{11} =-1,$ $u= x_1 ^2 = x_1^2 x_2=0.$ Therefore,  $\mu (u) =\{ x_1\}$  and $u$ is connected. Meantime, $\mu (u) =\{ x_1, x_2\}$  and $u$ is disconnected.

Throughout,  $\mathbb Z =: \{x \mid  x \hbox { is an integer}\}.$ $\mathbb R =: \{ x \mid x \hbox { is a real number}\}$.
$\mathbb N_0 =: \{x \mid  x \in \mathbb Z, x\ge 0\}.$
$\mathbb N =: \{x \mid  x \in \mathbb Z,  x>0\}$.  $F$ denotes the base field,   which is an algebraic closed field with  characteristic zero. $F^{*}=F\backslash\{0\}$. $\mathbb S_{n}$ denotes symmetric group, $n\in\mathbb N$. For any set $X$, $\mid X\mid$ is the cardinal of $X$. ${ {\rm int} }(a)$ means the biggest integer not greater than $a\in \mathbb R$.

\section {Relationships between  pure generalized Dynkin graph $\Gamma(V)$  and $\mathfrak L(V) $}

Obviously,  every pure generalized Dynkin graph is a graph in terms of  the graph theory (see \cite {Ha69}). Conversely,  assume that  $\Gamma$ is a  graph with vertex $\{1,  2,  \cdots,  n\}$. We define a matrix  $(p_{ij})_{n \times n}$ such that
$p_{ij}p_{ji} \neq 1$ (e.g. $p_{ij }=2$ and $p_{11} =p_{22} =\cdots = p_{nn} =2$) if and only if there exists an edge between $i$ and $j.$ Let $V$ be the braided vector space with braiding matrix $(p_{ij})_{n \times n}$. Then $\Gamma$ is a
pure generalized Dynkin graph of $V.$

By \cite [ Corollary  {2.5} and {5.2} ] {WWZZ17}  we have
\begin {Theorem}  \label{6.1}

{\rm (i)} $\Gamma (u)$ is a connected component of $\Gamma (V)$ with a non-zero monomial $u$ if and only if $\mu (u)$ is a maximal element in $\{\mu (v) \mid v \hbox { is a non-zero monomial in }  \mathfrak L(V) \}$ under order $\subseteq$.

{\rm (ii)} The following conditions are equivalent.

{\rm (a)}. $\Gamma (V)$ is connected

{\rm (b)}. $x_nx_{n-1} \cdots x_1  \in \mathfrak{L}(V)$.

 {\rm (c)}.  $x_1x_{2} \cdots x_n  \in \mathfrak{L}(V)$

{\rm (d)}.  there exists a non-zero monomial $u \in \mathfrak{L}(V)$ with $\mid \mu (u) \mid = \dim V$.

\end {Theorem}

\section {Relationships between augmented Dynkin  graph $\Gamma_a(V)$  and Nichols  Lie algebra $\mathfrak L^ - (V)$}

\begin {Proposition} \label {7.5} Assume that $h_i \in \{x_1,  \cdots,  x_n\}$ for $1\le i \le m$ and $u = h_1h_2\cdots h_m.$

{\rm (i)} If it is disconnected between monomial $u$ and monomial $v$ in $\Gamma_a (V)$ (i.e. ${p} _{x_i,  x_j} =1$  for any  $x_i \in \mu (u)$,  $x_j \in \mu (v)$ and $i\not= j$),  then $[u,  v]^- =0.$

{\rm (ii)} If $\Gamma_a (u)$  is disconnected in $\Gamma_a(V)$ or $ \mid \mu (u) \mid =1$  ),  then $\sigma (h_{1},  h_{2},  \cdots,   h_{m})=0$ for any  method $\sigma$ of adding bracket $[\ ]^-$ on $h_1,  h_2,  \cdots,  h_m.$

\end {Proposition}

\noindent {\it Proof.} {\rm (i)}  $u$ and $v$ are  commutative
(i.e. $uv= vu$ ) since $x_i x_j = x_jx_i$ for any  $x_i \in \mu (u)$,  $x_j \in \mu (v)$.

{\rm (ii)}
We show this by induction on $m$. $ [h_1,  h_2]^- =0$ for $m=2.$  For $m>2,  $ $\sigma (h_{1},  h_{2},  \cdots,  h_{m}) $

\noindent $= [\sigma _1 ( h_1 h_2 \cdots h_t ),  \sigma _2 ( h_{t+1} h_{t+2} \cdots h_m ) ]^-$.
If both $ h _1h_2\cdots h_t $ and $ h_{t+1} h_{t+2} \cdots h_m $ are connected  in $\Gamma_a (V)$,  then it is disconnected between  $h_1 h_2 \cdots h_t $ and $ h_{t+1} h_{t+2} \cdots h_m $  in $\Gamma_a (V)$. By Part {\rm (i)},   $\sigma (h_{1},  h_{2},  \cdots,  h_{m})=0$.
 If either  $h_1h_2 \cdots h_t $ or $ h_{t+1} h_{t+2} \cdots h_m$ is  disconnected  in $\Gamma_a (V)$,
then either  $\sigma _1 (h_1h_2\cdots h_t )=0$ or $ \sigma _2 (h_{t+1} h_{t+2} \cdots h_m)=0 $ by inductive hypothesis.

\hfill $\Box$


\begin{thebibliography}{BD99}



%\bibitem [AS02]{AS02} N. Andruskiewitsch,    H.-J. Schneider,  Pointed Hopf algebras. In New Directions in Hopf
%Algebras,  vol.43 ser MSRI Publications. Cambridge University Press(2002).


%\bibitem[ARS95] {ARS95}M. Auslander,  I. Reiten and S.O. Smal$\phi$,  Representation
%theory of Artin algebras,  Cambridge University Press,  1995.


\bibitem [AS10]{AS10} N. Andruskiewitsch,    H.-J. Schneider,     On the classification of finite-dimensional pointed Hopf algebras,
 Ann. Math.   {\bf 171}  (2010),    375-417.

\bibitem [AHS08] {AHS08} N. Andruskiewitsch,       I. Heckenberger and   H.J. Schneider,
  The Nichols algebra of a semisimple Yetter-Drinfeld module,
 Amer. J. Math. {\bf 132}  (2010),      1493-1547.

%\bibitem[An11] {An11} I. Angiono,        Nichols algebras of unidentified diagonal type,     arXiv:1108.5157.
%.SG); Quantum Algebra (math.QA)

\bibitem [AAB16] {AAB16} N. Andruskiewitsch, I Angiono, F. R. Bertone,
A finite-dimensional Lie algebra arising from a Nichols algebra of diagonal type (rank 2), arXiv:1603.09387.


\bibitem [BFM96]{BFM96} Y. Bahturin,  D. Fishman and S. Montgomery,  On the generalized
Lie structure of associative algebras,  J. Alg. {\bf 96} (1996),  27-48.

\bibitem [BFM01]{BFM01} Y. Bahturin,  D. Fischman and  S. Montgomery,
Bicharacter,  twistings and Scheunert's theorem for Hopf algebra,
J. Alg. {\bf 236} (2001),  246-276.

\bibitem [BMZP92] {BMZP92} Y. Bahturin,  D. Mikhalev,  M. Zaicev and V. Petrogradsky,
Infinite dimensional Lie superalgebras,  Walter de Gruyter Publ.
Berlin,  New York,  1992.


%\bibitem [GM03]{GM03} X. Gomez and S. Majid,   Braided Lie algebras and bicovariant
%differential calculi over coquasitriangular Hopf algebras,
%J. Alg. {\bf 261}(2003),  334--388.

\bibitem [GRR95]{GRR95} D. Gurevich,  A. Radul and V. Rubtsov,
Noncommutative differential geometry related to the Yang-Baxter
equation,  % Zap. Nauchn. Sem. S.-Peterburg Otdel. Mat. Inst. Steklov. (POMI) {\bf 199 } (1992); translation in
J. Math. Sci. {\bf 77 } (1995),  3051-3062.

\bibitem[Gu86] {Gu86} D. I. Gurevich,  The Yang-Baxter equation and the
generalization of formal Lie theory,  Dokl. Akad. Nauk SSSR,  {\bf
288} (1986),  797-801.
\bibitem[He06a]{He06a} I. Heckenberger,   Classification of arithmetic
root systems,   Adv. Math.  {\bf 220} (2009),   59-124.

\bibitem[He06b]{He06b} I. Heckenberger,   The Weyl-Brandt groupoid of a Nichols algebra
of diagonal type,   Invent. Math. {\bf 164} (2006),   175-188.

\bibitem[He05]{He05} I. Heckenberger,       Nichols algebras of diagonal type and arithmetic root systems,     Habilitation,   2005.

\bibitem[Ha69]{Ha69}  Frank Harary,
Graph Theory,  Addison¨CWesley,  USA,  1969

%\bibitem[HS08]{HS08} I. Heckenberger,   H.-J. Schneider,    {Root systems and Weyl groupoids for  Nichols algebras},   {arXiv:0807.0691}.

\bibitem[Hu78]{Hu78}  Humphreys,  James E. Introduction to Lie Algebras and Representation Theory,  Second printing,  revised. Graduate Texts in Mathematics,  9. Springer-Verlag,  New York,  1978.



\bibitem [Ka77]{Ka77} V. G. Kac,  Lie Superalgebras,   Adv. Math. {\bf 26} (1977),
8-96.

\bibitem   [Kh99] {Kh99} V. K. Kharchenko,    A Quantum analog of the
 poincar$\acute{e}$-Birkhoff-Witt theorem,    Algebra and
 Logic,    {\bf 38} (1999),    259-276


\bibitem   [LR95] {LR95} P. Lalconde and A. Ram,   Standard Lyndon bases of  Lie algebra and enveloping algebras, Trans.  AMS  {\bf 347} (1995), 1821-1830.



% \bibitem [Kh99b]{Kh99b} V. K. Kharchenko,
%An existence condition for multilinear quantum operations,  J. Alg.
%{\bf 217} (1999),  188--228.

\bibitem[Ma94b] {Ma94b}S. Majid,  Quantum and braided Lie algebras,  J. Geom. Phys.
{\bf 13} (1994),  307-356.

\bibitem[Ma95] {Ma95} S. Majid,  Foundations of quantum group,   Cambradge University
Press,  1995.

\bibitem [Pa98]{Pa98}  B. Pareigis,  On Lie algebras in the category of
Yetter-Drinfeld modules. Appl. Categ. Structures,   {\bf 6}
(1998),  151-175.

\bibitem [Sc79]{Sc79} M. Scheunert,  Generalized  Lie algebras,  J. Math. Phys.
{\bf 20} (1979), 712-720.


%\bibitem [Wo87]{Wo87} Woronowicz,  S.L.,  Compact matrix pseudogroups,  Commun. Math.
%Phys. 111 (1987),  613-665.


\bibitem[WZZ15a] {WZZ15a} W. Wu,    S. Zhang and   Y.-Z. Zhang,     Relationship between Nichols braided Lie
algebras and Nichols algebras,   J. Lie Theory {\bf 25} (2015),      45-63.

\bibitem[WZZ15b] {WZZ15b} W. Wu,    S. Zhang and   Y.-Z. Zhang,     On  Nichols (braided) Lie algebras,   Int. J. Math. {\bf 26} (2015),   1550082. % Also in arXiv:1409.3769.



\bibitem[WZZ16] {WZZ16} W. Wu,    S. Zhang and   Y.-Z. Zhang,
Finiteness of Nichols algebras and Nichols (braided) Lie algebras,
arXiv:1607.07955.


%\bibitem [Zh93]{Zh93} S.C. Zhang,  The Baer radical of generalized matrix rings,
%in Proc. of the Sixth SIAM Conf. on Parallel Processing for
%Scientific Computing,  pp.546--551,  Norfolk,  Virginia,  1993.
%Eds: R.F. Sincovec,  D.E. Keyes,  M.R. Leuze,  L.R. Petzold,  D.A.  Reed.

%\bibitem[ZZH03] {ZZH03} S.C. Zhang,  Y.Z. Zhang and  Y.Y. Han,
%Duality theorems for infinite braided Hopf algebras,  math.QA/0309007.

\bibitem [ZZ03]{ZZ03} S. C. Zhang and  Y.-Z. Zhang,
Braided m-Lie algebras,   Lett. Math. Phys.  {\bf 70} (2004),  155-167.


% \bibitem [ZWTZ]{ZWTZ} S. Zhang,   W. Wu,   Zhengtang Tan and   Y.-Z. Zhang,
%Nichols algebras over classical Weyl groups,   Fomin-Kirillov algebras and Lyndon basis,
%     arXiv:1307.8227

\bibitem [WWZZ17]{WWZZ17} Weicai Wu, Jing Wang, Shouchuan Zhang, Yao-Zhong Zhang,
Structures of Nichols (braided) Lie algebras of diagonal type, Journal of Lie Theory,
 28 (2018), 357-380. Also see  arXiv:1704.06810.

\end{thebibliography}
\end {document}